\documentclass[12pt,twoside]{amsart}
\usepackage{amsfonts}
\usepackage{amssymb}
\usepackage{amsmath}
\usepackage{fancyhdr}

\theoremstyle{plain}
\newtheorem{thm}{Theorem}[section]
\newtheorem{theorem}[thm]{Theorem}
\newtheorem{lemma}[thm]{Lemma}
\newtheorem{corollary}[thm]{Corollary}
\newtheorem{proposition}[thm]{Proposition}

\theoremstyle{definition}

\newtheorem{defn-thm}[thm]{Definition-Theorem}

\numberwithin{equation}{section}
\catcode`\@=11
\def\opn#1#2{\def#1{\mathop{\kern0pt\fam0#2}\nolimits}}
\def\underrightarrow{\mathpalette\underrightarrow@}
\def\underrightarrow@#1#2{\vtop{\ialign{$##$\cr
 \hfil#1#2\hfil\cr\noalign{\nointerlineskip}%
 #1{-}\mkern-6mu\cleaders\hbox{$#1\mkern-2mu{-}\mkern-2mu$}\hfill
 \mkern-6mu{\to}\cr}}}

\def\underleftarrow{\mathpalette\underleftarrow@}
\def\underleftarrow@#1#2{\vtop{\ialign{$##$\cr
 \hfil#1#2\hfil\cr\noalign{\nointerlineskip}#1{\leftarrow}\mkern-6mu
 \cleaders\hbox{$#1\mkern-2mu{-}\mkern-2mu$}\hfill
 \mkern-6mu{-}\cr}}}
\let\amp@rs@nd@\relax
\newdimen\ex@
\newdimen\bigaw@
\newdimen\minaw@
\newdimen\minCDaw@
\newif\ifCD@
\def\minCDarrowwidth#1{\minCDaw@#1}

\def\@CD{\def\A##1A##2A{\llap{$\vcenter{\hbox
 {$\scriptstyle##1$}}$}\Big\uparrow\rlap{$\vcenter{\hbox{%
$\scriptstyle##2$}}$}&&}%
\def\V##1V##2V{\llap{$\vcenter{\hbox
 {$\scriptstyle##1$}}$}\Big\downarrow\rlap{$\vcenter{\hbox{%
$\scriptstyle##2$}}$}&&}%
\def\={&\hskip.5em\mathrel
 {\vbox{\hrule width\minCDaw@\vskip3\ex@\hrule width
 \minCDaw@}}\hskip.5em&}%
\def\verteq{\Big\Vert&&}%
\def\noarr{&&}%
\def\vspace##1{\noalign{\vskip##1\relax}}\relax\let\amp@rs@nd@&\iffalse}\fi
 \CD@true\vcenter\bgroup\relax\let\\=\cr\iffalse}\fi\tabskip\z@skip\baselineskip20\ex@
 &\hfill$\m@th##$\hfill\cr}
\def\@endCD{\cr\egroup\egroup}
\def\>#1>#2>{\amp@rs@nd@\setbox\z@\hbox{$\scriptstyle
 \;{#1}\;\;$}\setbox\@ne\hbox{$\scriptstyle\;{#2}\;\;$}\setbox\tw@
 \hbox{$#2$}\ifCD@
 \global\bigaw@\minCDaw@\else\global\bigaw@\minaw@\fi
 \ifdim\wd\z@>\bigaw@\global\bigaw@\wd\z@\fi
 \ifdim\wd\@ne>\bigaw@\global\bigaw@\wd\@ne\fi
 \ifCD@\hskip.5em\fi
 \ifdim\wd\tw@>\z@
 \mathrel{\mathop{\hbox to\bigaw@{\rightarrowfill}}\limits^{#1}_{#2}}\else
 \mathrel{\mathop{\hbox to\bigaw@{\rightarrowfill}}\limits^{#1}}\fi
 \ifCD@\hskip.5em\fi\amp@rs@nd@}
\def\<#1<#2<{\amp@rs@nd@\setbox\z@\hbox{$\scriptstyle
 \;\;{#1}\;$}\setbox\@ne\hbox{$\scriptstyle\;\;{#2}\;$}\setbox\tw@
 \hbox{$#2$}\ifCD@
 \global\bigaw@\minCDaw@\else\global\bigaw@\minaw@\fi
 \ifdim\wd\z@>\bigaw@\global\bigaw@\wd\z@\fi
 \ifdim\wd\@ne>\bigaw@\global\bigaw@\wd\@ne\fi
 \ifCD@\hskip.5em\fi
 \ifdim\wd\tw@>\z@
 \mathrel{\mathop{\hbox to\bigaw@{\leftarrowfill}}\limits^{#1}_{#2}}\else
 \mathrel{\mathop{\hbox to\bigaw@{\leftarrowfill}}\limits^{#1}}\fi
 \ifCD@\hskip.5em\fi\amp@rs@nd@}

\def\@CDS{\def\A##1A##2A{\llap{$\vcenter{\hbox
 {$\scriptstyle##1$}}$}\Big\uparrow\rlap{$\vcenter{\hbox{%
$\scriptstyle##2$}}$}&}%
\def\V##1V##2V{\llap{$\vcenter{\hbox
 {$\scriptstyle##1$}}$}\Big\downarrow\rlap{$\vcenter{\hbox{%
$\scriptstyle##2$}}$}&}%
\def\={&\hskip.5em\mathrel
 {\vbox{\hrule width\minCDaw@\vskip3\ex@\hrule width
 \minCDaw@}}\hskip.5em&}
\def\verteq{\Big\Vert&}
\def\novarr{&}
\def\noharr{&&}
\def\SE##1E##2E{\slantedarrow(0,18)(4,-3){##1}{##2}&}
\def\SW##1W##2W{\slantedarrow(24,18)(-4,-3){##1}{##2}&}
\def\NE##1E##2E{\slantedarrow(0,0)(4,3){##1}{##2}&}
\def\NW##1W##2W{\slantedarrow(24,0)(-4,3){##1}{##2}&}
\def\slantedarrow(##1)(##2)##3##4{%
\thinlines\unitlength1pt\lower 6.5pt\hbox{\begin{picture}(24,18)%
\put(##1){\vector(##2){24}}%
\put(0,8){$\scriptstyle##3$}%
\put(20,8){$\scriptstyle##4$}%
\end{picture}}}
\def\vspace##1{\noalign{\vskip##1\relax}}\relax\let\amp@rs@nd@&\iffalse}\fi
 \CD@true\vcenter\bgroup\relax\let\\=\cr\iffalse}\fi\tabskip\z@skip\baselineskip20\ex@
 &\hfill$\m@th##$\hfill\cr}
\def\@endCDS{\cr\egroup\egroup}
\newdimen\TriCDarrw@
\newif\ifTriV@

\def\@TriCDV{\TriV@true\def\TriCDpos@{6}\@TriCD}
\def\@TriCDA{\TriV@false\def\TriCDpos@{10}\@TriCD}
\def\@TriCD#1#2#3#4#5#6{%
\setbox0\hbox{$\ifTriV@#6\else#1\fi$} \TriCDarrw@=\wd0
\advance\TriCDarrw@ 24pt \advance\TriCDarrw@ -1em
\def\SE##1E##2E{\slantedarrow(0,18)(2,-3){##1}{##2}&}
\def\SW##1W##2W{\slantedarrow(12,18)(-2,-3){##1}{##2}&}
\def\NE##1E##2E{\slantedarrow(0,0)(2,3){##1}{##2}&}
\def\NW##1W##2W{\slantedarrow(12,0)(-2,3){##1}{##2}&}
\def\slantedarrow(##1)(##2)##3##4{\thinlines\unitlength1pt
\lower 6.5pt\hbox{\begin{picture}(12,18)%
\put(##1){\vector(##2){12}}%
\put(-4,\TriCDpos@){$\scriptstyle##3$}%
\put(12,\TriCDpos@){$\scriptstyle##4$}%
\end{picture}}}
\def\={\mathrel {\vbox{\hrule
   width\TriCDarrw@\vskip3\ex@\hrule width
   \TriCDarrw@}}}
\def\>##1>>{\setbox\z@\hbox{$\scriptstyle
 \;{##1}\;\;$}\global\bigaw@\TriCDarrw@
 \ifdim\wd\z@>\bigaw@\global\bigaw@\wd\z@\fi
 \hskip.5em
 \mathrel{\mathop{\hbox to \TriCDarrw@
{\rightarrowfill}}\limits^{##1}}
 \hskip.5em}
\def\<##1<<{\setbox\z@\hbox{$\scriptstyle
 \;{##1}\;\;$}\global\bigaw@\TriCDarrw@
 \ifdim\wd\z@>\bigaw@\global\bigaw@\wd\z@\fi
 \mathrel{\mathop{\hbox to\bigaw@{\leftarrowfill}}\limits^{##1}}
 }
 \CD@true\vcenter\bgroup\relax\let\\=\cr\iffalse}\fi
 \tabskip\z@skip\baselineskip20\ex@
 \ifTriV@
 \halign\bgroup
 &\hfill$\m@th##$\hfill\cr
#1&\multispan3\hfill$#2$\hfill&#3\\
&#4&#5\\
&&#6\cr\egroup%
\else
 \halign\bgroup
 &\hfill$\m@th##$\hfill\cr
&&#1\\%
&#2&#3\\
#4&\multispan3\hfill$#5$\hfill&#6\cr\egroup \fi}
\def\@endTriCD{\egroup}

\newcounter{Myenumi}
{\begin{list}{}{\usecounter{Myenumi}%
\settowidth{\leftmargin}{2.n}\settowidth{\labelwidth}{2.n}%
\setlength{\labelsep}{0pt}}}{\end{list}}
\newcounter{Myenumii}
{\begin{list}{}{\usecounter{Myenumii}%
\settowidth{\leftmargin}{a)n}\settowidth{\labelwidth}{a)n}%
\setlength{\labelsep}{0pt}}}{\end{list}}
\newcounter{Myenumiii}
{\begin{list}{}{\usecounter{Myenumiii}%
\settowidth{\leftmargin}{iv.n}\settowidth{\labelwidth}{iv.n}%
\setlength{\labelsep}{0pt}}}{\end{list}}

{\end{list}}

{\begin{list}{}{%
\settowidth{\leftmargin}{2.n}\settowidth{\labelwidth}{2.n}%
\setlength{\labelsep}{0pt}}}{\end{list}}

\newsymbol\onto 1310

\begin{document}

\title{Regularity of Leray-Hopf solutions to Navier-Stokes equations
(I)-Critical  regularity in weak spaces}

\author{Jian Zhai \\ \tiny{Department of Mathematics, Zhejiang University, Hangzhou
310027, PRC} }
\maketitle

\footnotetext{This work is supported by NSFC No.10571157. email:
jzhai@zju.edu.cn}




\begin{abstract}
We consider the  regularity of  Leray-Hopf solutions to impressible
Navier-Stokes equations on critical case $u\in
L^2_w(0,T;L^\infty(\mathbb R^3))$. By a new embedding inequality in
Lorentz space we prove that if $\|u\|_{L^2_w(0,T;L^\infty(\mathbb
R^3))}$ is small then as a Leray-Hopf solution $u$ is regular.
Particularly, an open problem proposed in \cite{[KK]} is solved.

\end{abstract}

\section{Introduction}

We consider the regularity of weak solutions to impressible
Navier-Stokes equations
\begin{equation} \label{1.1}
\left\{ \begin{aligned}
\partial_tu-\Delta u+u\cdot\nabla u+\nabla p=0,\quad\text{in}\quad \mathbb R^3\times(0,T) \\
\text{div}u=0,\quad\text{in}\quad \mathbb R^3\times(0,T)
\end{aligned} \right.
\end{equation}
where  $u$ and $p$ denote the unknown velocity and pressure of
incompressible fluid respectively.  $u: \,\, (x,t)\in\mathbb
R^3\times(0,T)\to\mathbb R^3$ is called a weak solution of
(\ref{1.1}) if it is a Leray-Hopf solution. Precisely, it
satisfies
\begin{eqnarray*}
(1)&& \qquad u\in L^\infty(0,T;L^2(\mathbb R^3))\cap L^2(0,T;H^{1}(\mathbb R^3)),\\
(2)&& \qquad \text{div} u=0 \quad\text{in}\quad \mathbb R^3\times(0,T),\\
(3)&& \qquad \int_0^T\int_{\mathbb
R^3}\{-u\cdot\partial_t\phi+\nabla u\cdot\nabla\phi+(u\cdot\nabla
u)\cdot\phi\}dxdt=0
\end{eqnarray*}
for all $\phi\in C^\infty_0(\mathbb R^3\times(0,T))$ with div$\phi=0$ in $\mathbb R^3\times(0,T)$.\\

In this paper, we prove the following critical regularity of the
Leray-Hopf solutions to the Navier-Stokes equations in weak spaces,
which was an open
problem proposed in \cite{[KK]}.\\

\begin{theorem} \label{thm1.1}
There is a constant $\epsilon>0$ such that if $u$ is a weak solution
of the Navier-Stokes equations (\ref{1.1}) in $\mathbb
R^3\times(0,T)$ and if
\begin{eqnarray*}
\|u\|_{L^2_w(0,T;L^\infty(\mathbb R^3))}\leq\epsilon
\end{eqnarray*}
then $u$ is regular in $\mathbb R^3\times(0,T]$.
\end{theorem}\vspace{0.5cm}

Here $L^p_w(0,T;L^q(\Omega))$ ($1<p<\infty$, $1\leq q\leq \infty$)
denote the spaces of functions $v:\,\, (x,t)\in\Omega\times(0,T)\to
\mathbb R^3$ with
$$
\|v\|_{L^p_w(0,T;L^q(\Omega))}:=\sup_\sigma \sigma
|\{t\in(0,T):\quad
\|v(\cdot,t)\|_{L^q(\Omega)}>\sigma\}|^{1/p}<\infty.
$$
It is known that the weak spaces $L^p_w$ are special cases of the
more general Lorentz spaces $L^{p,r}$ and $L^p_w=L^{p,\infty}$ (see
 \cite{[BL]}).\\

 As a corollary of  Theorem \ref{thm1.1}, we have that there is
a constant $C>0$ such that if $u$ is a weak solution of (\ref{1.1})
and
$$
|u(x,t)|\leq \frac C{(T-t)^{1/2}},\quad\forall (x,t)\in \mathbb
R^3\times (T-R,T)
$$
then $u$ is bounded in $\mathbb R^3\times (T-R,T]$.\\

 Combining our  Theorem \ref{thm1.1} with the former results of Sohr \cite{[So]}, Kim and
Kozono \cite{[KK]}, we have

\begin{corollary} \label{cor1.2}
 For all $r\in[3,\infty]$, there is a constant
$\epsilon>0$ depending only on $r$, such that if $u$ is a weak
solution of the Navier-Stokes equations (\ref{1.1}) and if
$$
\|u\|_{L^s_w(0,T;L^r_w(\mathbb R^3))}\leq\epsilon,\quad
\text{with}\quad \frac2s+\frac3r=1,
$$
then $u$ is regular in $\mathbb R^3\times(0,T]$.
\end{corollary}\vspace{0.5cm}

Since Leray(1934)\cite{[L]} and Hopf(1951)\cite{[H]} proved the
global existence of weak solutions,  it has been a fundamental open
problem to prove the uniqueness and regularity of weak solutions to
the Navier-Stokes equations. For $3<r<\infty$, Corollary
\ref{cor1.2} in $L^s_w(0,T;L^r_w(\mathbb R^3))$ were proved by Sohr
\cite{[So]}. Corollary \ref{cor1.2} in $r=3$ was proved by Kim and
Kozono \cite{[KK]}. On the other hand, similar results in Lebesgue
spaces on $\Omega\times(0,T)$ have been proved by Serrin \cite{[S]},
Struwe \cite{[St]} and Takahashi \cite{[T]}, and similar results in
Lebesgue spaces on $\mathbb R^3\times(0,T)$ have been proved by Giga
\cite{[G]}, E.B. Fabes, B.F. Jones, N.M. Rivere \cite{[FJR]},
Kozono, Taniuchi \cite{[KT]} and Iskauriaza, Ser\"{e}gin, Shverak
\cite{[ISS]} ( see also W. von Wahl \cite{[W]}).\\

Notice that the global case of the open problem proposed by Kim and
Kozono in [KK pp.87 line 12-14] is solved by using Theorem
\ref{thm1.1}. The local case of the open problem was claimed in [T].
But as pointed by Kim and Kozono in [KK pp.99 line 9-11], the
critical local case can
not be treated by the method given in [T] and developed in [KK].\\

To prove Theorem \ref{thm1.1}, a key step is to prove a priori
estimate for vorticity equation (see Proposition \ref{pro3.1}),
where  we estimate the nonlinear terms by $\|v\|_Q$ and the norm of
$u$ in Lorentz space. To this aim, we first prove a new embedding
inequality in
Lorentz space in section 2.\\

\section{ Embedding inequality in Lorentz space }

Let $\chi\in C_0^\infty(B_{4/3}(0))$ and $\varphi\in
C_0^\infty(B_{8/3}(0)\setminus B_{3/4}(0))$ be the Littlewood-Paley
dyadic decomposition that satisfy (see [C]):
\begin{equation}\label{2.1}
\chi(\xi)+\sum_{q\geq0}\varphi(2^{-q}\xi)=1,\quad
\frac13\leq\chi^2(\xi)+\sum_{q\geq0}\varphi^2(2^{-q}\xi)\leq1,\quad\forall
\xi\in\mathbb R^3.
\end{equation}
Denote
$$
\Delta_{-1}v=\mathcal{F}^{-1}[\chi(\xi)\mathcal{F}[v](\xi)],\quad
\Delta_qv=\mathcal{F}^{-1}[\varphi(2^{-q}\xi)\mathcal{F}[v](\xi)],\quad\forall
q\geq0
$$
and  define
\begin{equation}
\begin{aligned}\label{2.2}
V(Q[0,T]):=\{v\in L^2(0,T;H^1(\mathbb R^3)):\quad \|v\|_{Q[0,T]}<\infty\}\\
\|v\|_{Q[0,T]}^2=\sum_{q\geq -1}\sup_{0\leq t<
T}\frac12\int_{\mathbb
R^3}|\Delta_qv(x,t)|^2dx+\int_0^T\int_{\mathbb R^3}|\nabla
v(x,t)|^2dxdt.
\end{aligned}
\end{equation}
We shall use the notation $\|v\|_Q$ and $V(Q)$ to denote
$\|v\|_{Q[0,1]}$ and $V(Q[0,1])$ respectively.

\begin{lemma}\label{lem2.1}
There is a constant $C>0$, such that for all $f\in
L^{2,\infty}(0,1)$ and $v\in V(Q)$,
\begin{equation}\label{2.3}
\sum_{q\geq-1}\sum_{q-2\leq j\leq q+4}\int_0^1|f(t)|\int_{\mathbb
R^3}|\Delta_jv(x,t)||\nabla\Delta_q v(x,t)|dxdt\leq
C\|f\|_{L^{2,\infty}}\|v\|_Q^2.
\end{equation}

\end{lemma}\vspace{0.5cm}

$Proof.$ Step 1. Note that the weak space $L^2_w(0,1)$ is
equivalent to the Lorentz space $L^{2,\infty}(0,1)$, and the norm
on $L^{2,\infty}(0,1)$ can be defined equivalently by
$$
\|f\|_{L^{2,\infty}(0,1)}=\sup\{|E|^{-1/2}\int_E|f(t)|dt;\quad E\in
\mathcal L\}
$$
where $\mathcal L$ is the collection of all Lebesgue measurable sets
with a positive measure (see \cite{[M]}). Instead of the Lebesgue
measurable sets, the original version in \cite{[M]}(18.5) used the
collection of all Borel sets with a positive measure. Since for all
Lebesgue measurable sets $E$ and $1<p<\infty$
\begin{eqnarray*}
&& \int_E|f(t)|dt=\int_{E\cap\{|f|\geq
\sigma\}}|f(t)|dt+\int_{E\setminus\{|f|\geq \sigma\}}|f(t)|dt\\
&& \leq \int_\sigma^\infty \lambda^p|\{t\in E:\,\,
|f(t)|>\lambda\}|\frac{d\lambda}{\lambda^p}+\sigma |E|\\
&& \leq C\|f\|_{L^p_w(E)}|E|^{1/p'}
\end{eqnarray*}
by taking $\sigma=\|f\|_{L^p_w(E)}|E|^{\frac{1}{(1-p)p'}}$, nothing is lost when we use Lebesgue measurable sets to replace Borel measurable sets.\\

It is known that $L^{2,\infty}(0,1)$ is the dual space of
$L^{2,1}(0,1)$, where for $g\in L^{2,1}(0,1)$ the norm is defined by
the infimum of $\sum_{j\geq0}|c_j|$, the sums of the coefficients of
the atom decomposition
$$
g(t)=\sum_{j\geq0}c_ja_j(t)
$$
over all possible expansions of $g$.\\

\vspace{0.5cm}

Step 2. Note that for $q\geq -1$
\begin{equation}\label{2.4}
\begin{aligned}
\|\nabla\Delta_q v\|_{L^2(\mathbb R^3)}=(\int_{\mathbb R^3}\sum_{1\leq j\leq3}|\nabla_{x_j}\Delta_qv|^2dx)^{1/2} \\
=(\int_{\mathbb R^3}\sum_{1\leq j\leq3}|i\xi_j\mathcal{F}[\Delta_qv](\xi)|^2d\xi)^{1/2}\\
\leq (\frac{8}3)2^q\|\Delta_qv\|_{L^2(\mathbb R^3)},
\end{aligned}
\end{equation}
and for $q\geq 0$
\begin{equation}\label{2.5}
\|\nabla\Delta_q v\|_{L^2(\mathbb R^3)}\geq (\frac34)
2^{q}\|\Delta_q v\|_{L^2(\mathbb R^3)}.
\end{equation}

Denote
$$
M(v)=\sup_{0\leq t<1}( \|\Delta_qv(t)\|_{L^2(\mathbb
R^3)}\|\Delta_jv(t)\|_{L^2(\mathbb R^3)})
$$
 and for
$k=1,2,3,...,$ define
$$
E_k=\{t\in (0,1): 2^{-k}< (M(v))^{-1}(
\|\Delta_qv(t)\|_{L^2(\mathbb R^3)}\|\Delta_jv(t)\|_{L^2(\mathbb
R^3)})\leq 2^{-(k-1)}\}.
$$
Since $\|v\|_Q$ is bounded, $M(v)$ is bounded and $E_k$ are Lebesgue
measurable.

Note that for $t\in E_k$
\begin{equation}\label{2.6}
\begin{aligned}
 ( \|\Delta_qv(t)\|_{L^2(\mathbb
R^3)}\|\Delta_jv(t)\|_{L^2(\mathbb
R^3)})\leq 2^{-(k-1)}M(v)\\
< \frac2{|E_k|}\int_{E_k}( \|\Delta_qv(t)\|_{L^2(\mathbb
R^3)}\|\Delta_jv(t)\|_{L^2(\mathbb R^3)})dt.
\end{aligned}
\end{equation}

\vspace{0.5cm}
 Step 3. Denote
$$
h(t)=( \|\Delta_qv(t)\|_{L^2(\mathbb
R^3)}\|\Delta_jv(t)\|_{L^2(\mathbb R^3)})
$$
and notice that
\begin{eqnarray}\label{2.7}
&& \int_0^1|f(t)|\int_{\mathbb R^3}|\Delta_jv(x,t)||\nabla \Delta_qv(x,t)|dxdt \nonumber\\
&& \leq\int_0^1|f(t)|(\int_{\mathbb R^3}|\Delta_jv(x,t)|^2dx)^{1/2}(\int_{\mathbb R^3}|\nabla \Delta_qv(x,t)|^2dx)^{1/2}dt \nonumber\\
&& \leq(\frac83)2^q\int_0^1|f(t)|h(t)dt\quad \text{(by (\ref{2.4}))} \nonumber\\
&& \leq (\frac83)2^{q+1}\sum_{k\geq 1}|E_k|^{-1}\int_{E_k}|f(t)|dt\int_{E_k}h(t)dt\quad \text{(by (\ref{2.6}))} \nonumber\\
&& \leq (\frac83)2^{q+1}\|f\|_{L^{2,\infty}}  \sum_{k\geq 1}\frac1{|E_k|^{1/2}}\int_{E_k}h(t)dt\quad \text{(by step 1)}\nonumber\\
&& \leq(\frac83)2^{q+1}\|f\|_{L^{2,\infty}}  \sum_{k\geq 1}(\sup_{E_k}h)^{1/2}(\int_{E_k}h(t)dt)^{1/2}\nonumber\\
&& \leq (\frac83)2^{q+1}\|f\|_{L^{2,\infty}}  (\sum_{k\geq 1}\sup_{E_k}h)^{1/2}(\sum_{k\geq 1}\int_{E_k}h(t)dt)^{1/2}\nonumber\\
&& \leq  (\frac83)2^{q+1}\|f\|_{L^{2,\infty}}\sqrt{2}
M(v)^{1/2}(\int_0^1h(t)dt)^{1/2}. \quad \text{(by (\ref{2.6}))}
\end{eqnarray}
For $j,\,\,q\geq0$, by (\ref{2.5}) we have
\begin{eqnarray*}
&& \text{ the right of (\ref{2.7})}\\
&& \leq C\|f\|_{L^{2,\infty}} M(v)^{1/2}(\int_0^1(\int_{\mathbb
R^3}|\nabla\Delta_jv(x,t)|^2dx)^{1/2}(\int_{\mathbb
R^3}|\nabla\Delta_qv(x,t)|^2dx)^{1/2}dt)^{1/2},
\end{eqnarray*}
and
for $j=-1$, $q\geq 0$
\begin{eqnarray*}
&& \text{the right of (\ref{2.7})}\\
&& \leq C\|f\|_{L^{2,\infty}} M(v)^{1/2}(\int_0^1\int_{\mathbb
R^3}|\Delta_{-1}v(x,t)|^2dxdt)^{1/4}(\int_0^1\int_{\mathbb
R^3}|\nabla\Delta_{q}v(x,t)|^2dxdt)^{1/4}.
\end{eqnarray*}
So, by (\ref{2.1}) we have
\begin{equation}\label{2.8}
\begin{aligned}
\sum_{q\geq-1}\sum_{q-2\leq j\leq q+4}\int_0^1|f(t)|dt\int_{\mathbb R^3}|\Delta_jv(x,t)||\nabla \Delta_qv(x,t)|dx\\
\leq C\|f\|_{L^{2,\infty}}\sum_{q\geq0}(\sup_{0\leq t< 1}\|\Delta_qv(t)\|_{L^2(\mathbb R^3)})  (\int_0^1\int_{\mathbb R^3}|\nabla\Delta_qv(x,t)|^2dxdt)^{1/2}\\
+C\|f\|_{L^{2,\infty}}(\sup_{0\leq t< 1}\int_{\mathbb R^3}|\Delta_{-1}v(x,t)|^2dx)^{1/2}(\int_0^1\int_{\mathbb R^3}|\Delta_{-1}v(x,t)|^2dxdt)^{1/2}\\
\leq C\|f\|_{L^{2,\infty}}\{(\sum_{q\geq-1}\sup_{0\leq t< 1}\int_{\mathbb R^3}|\Delta_qv(x,t)|^2dx)^{1/2}(\int_0^1\int_{\mathbb R^3}|\nabla v(x,t)|^2dxdt)^{1/2}\\
+\sum_{q\geq-1}\sup_{0\leq t< 1}\int_{\mathbb
R^3}|\Delta_qv(x,t)|^2dx\}\\
\leq C\|f\|_{L^{2,\infty}}\|v\|_Q^2.\qed
\end{aligned}
\end{equation}

\section{ Proof of theorem \ref{thm1.1}}

Without loss generality, we assume $T=1$. We consider the Cauchy
problem for the vorticity equation which follows the Navier-Stokes
equations (\ref{1.1})
\begin{equation}\label{3.1}
\left\{\begin{aligned}
\partial_tv-\Delta v+\text{div}(Bv)=0,\quad\forall (x,t)\in\mathbb R^3\times(0,1)\\
v(x,0)=v_0(x),\quad\forall x\in\mathbb R^3,
\end{aligned}\right.
\end{equation}
where $Bv=v\otimes u-u\otimes v$, and $v=\text{curl}\,\,u$. The
following a priori estimate for (\ref{3.1}) will be proved in
section 4.

\begin{proposition}\label{pro3.1} There exists $\epsilon>0$ such that if
\begin{equation}\label{3.2}
\|u\|_{L^{2,\infty}(0,1;L^\infty(\mathbb R^3))}\leq\epsilon
\end{equation}
 then  for all $t_1\in (0,1]$, for all solutions $v$ of (\ref{3.1}) in $V(Q[0,t_1])$, we have
\begin{equation}\label{3.3}
\|v\|_{Q[0,t_1]}^2\leq C\|v_0\|_{L^2(\mathbb R^3)}^2
\end{equation}
where the constant $C$ is independent of $v$ and $t_1$.
\end{proposition}\vspace{0.5cm}

$Proof\,\,\,of\,\,\, Theorem \ref{1.1}$: Note that the weak space
$L^2_w(0,1)$ is equivalent to the Lorentz space $L^{2,\infty}(0,1)$.
Proposition \ref{pro3.1} implies a priori estimate for the solutions
of (\ref{3.1})  provided that
$$
\|u\|_{L^2_w(0,1;L^\infty(\mathbb R^3))}\leq\epsilon.
$$

If $u$ is a Leray-Hopf solution to (\ref{1.1}), then $u\in
L^2(0,1;H^{1}(\mathbb R^3))$. So for any $\delta_0>0$ there is
$\delta\in(0,\delta_0)$ such that $\|\nabla u(\delta)\|_{L^2(\mathbb
R^3)}<\infty$. Take
$$
v_0(x)=\text{curl}\,\,\, u(x,\delta)$$ and consider the Cauchy
problem
\begin{equation}\label{3.4}
\left\{\begin{aligned}
\partial_tv-\Delta v+\text{div}(Bv)=0,
\quad\text{in}\quad \mathbb
R^3\times(\delta,1)\\
v(x,\delta)=v_0(x),\quad\text{in}\quad\mathbb R^3.
\end{aligned}\right.
\end{equation}

The solution $v$ of (3.4) is regular at least in a short time
interval $(\delta,t_1)$ with $t_1\leq 1$. So for any small
$\delta_1>0$, $v\in V(Q[\delta,t_1-\delta_1])$. We can use the a
priori estimate in Proposition \ref{pro3.1} to get
$$
\|v\|_{Q[\delta,t_1]}^2\leq \limsup_{\delta_1\to
0+}\|v\|^2_{Q[\delta,t_1-\delta_1]}\leq C\|v_0\|_{L^2(\Bbb R^3)}^2
$$
because the constant $C$ is independent of $\delta_1$. Then for all
$t\in[\delta,t_1)$, $\|v(t)\|_{L^2(\mathbb R^3)}$ is uniformly
bounded. So $v$ is regular at $t=t_1$. Similarly by the initial data
$v(x,t_1)$ and so on we can prove that the solution $v$ of
(\ref{3.4}) is regular in $(\delta,1]$ provided that
$$
\|u\|_{L^2_w(0,1;L^\infty(\mathbb R^3))}\leq\epsilon.
$$

Let $\delta_0\to 0+$.  So $v(x,t)$ and $u(x,t)$ are regular for
$t\in (0,1]$ provided that
$$
\|u\|_{L^2_w(0,1;L^\infty(\mathbb R^3))}\leq\epsilon.
$$
Thus we proved Theorem \ref{thm1.1}.\qed

\section{Proof of proposition \ref{pro3.1}}

Without loss generality, we assume $t_1=1$. We introduce Bony's
paraproduct from the Littlewood-Paley analysis. We denote
$$
S_ju=\sum_{-1\leq k\leq j-1}\Delta_ku,\quad
\Delta_j(u)=S_{j+1}(u)-S_j(u).
$$
For the product $uv$ of $u$ and $v$, we shall decompose it as the
sum
$$
uv=T_uv+T_vu+R(u,v)
$$
 of  paraproducts
$$
T_uv:=\sum_{j\geq1}S_{j-1}u\Delta_jv,\quad
T_vu:=\sum_{j\geq1}\Delta_juS_{j-1}v,
$$
 and  remainder
 $$
 R(u,v):=\sum_{j\geq-1}\sum_{j-1\leq k\leq j+1}\Delta_ku\Delta_jv,
 $$
where
$$
S_jv=\sum_{-1\leq k\leq j-1}\Delta_kv,\quad
S_0v=\Delta_{-1}v,\quad S_{-1}v=0.
$$

Note that for $q\geq-1$,
$$
\Delta_q(T_vu)=\Delta_q(\sum_{q-2\leq j\leq q+4}\Delta_juS_{j-1}v)
$$
and
$$
\Delta_q(R(u,v))=\Delta_q(\sum_{j\geq q-3}\sum_{k=j-1}^{
j+1}\Delta_ku\Delta_jv),
$$
because (see \cite{[M]} Lemma 16), for example,
$$
\text{spt}(\Delta_juS_{j-1}v)\subset \{(\frac34-\frac23)2^j\leq
|\xi|\leq (\frac83+\frac23)2^j\},\quad \forall j\geq 2,
$$
and
$$
\text{spt}(\Delta_q)\subset \{(\frac34)2^q\leq |\xi|\leq
(\frac83)2^q\},\quad\forall q\geq 0,
$$
the necessary condition of $\text{spt}(\Delta_juS_{j-1}v)\cap
\text{spt}(\Delta_q)\neq\emptyset$ is $q-2\leq j\leq q+4$.

\vspace{0.5cm}

Step 1. Applying the operator $\Delta_q$ to (\ref{3.1}) we get
\begin{equation}\label{4.1}
\left\{\begin{aligned}
\partial_t\Delta_qv-\Delta \Delta_qv-\text{div} \Delta_q(u\otimes v-v\otimes u)=0,\quad\forall (x,t)\in\mathbb R^3\times(0,1)\\
\Delta_qv(x,0)=\Delta_qv_0(x),\quad\forall x\in\mathbb R^3.
\end{aligned}\right.
\end{equation}

Taking inner products with $\Delta_qv$ in the two sides of the
 equations (\ref{4.1}), we have
\begin{equation}\label{4.2}
\begin{aligned}
\frac12(\partial_t-\Delta)|\Delta_qv(x,t)|^2+|\nabla\Delta_qv(x,t)|^2\\
=-[\text{div}(\Delta_qR(B,v))]\cdot \Delta_qv(x,t)-[\text{div}(\Delta_qT_vB)]\cdot \Delta_qv(x,t)\\
-[\text{div}(\Delta_qT_Bv)]\cdot \Delta_qv(x,t),\quad\forall (x,t)\in\mathbb R^3\times(0,1),\\
|\Delta_qv(x,0)|^2=|\Delta_qv_0(x)|^2,\quad\forall x\in\mathbb
R^3.
\end{aligned}\end{equation}
Here the notations $R(B,v)$, $T_vB$, $T_Bv$ may be understood as
$R(u,v)$, $T_vu$, $T_uv$.

 Integrating (\ref{4.2}) over $\mathbb R^3\times[0,1]$ we have

\begin{equation}\label{4.3}
\begin{aligned}
\sup_{0\leq t< 1}\frac12\int_{\mathbb R^3}|\Delta_qv(x,t)|^2dx-\int_{\mathbb R^3}|\Delta_qv_0(x)|^2dx+\int_0^1dt\int_{\mathbb R^3}|\nabla\Delta_qv(x,t)|^2dx\\
\leq 2\int_0^1dt\int_{\mathbb R^3}\{\Delta_qR(B,v)+\Delta_qT_vB+\Delta_qT_Bv  \}\cdot\nabla\Delta_qv(x,t)dx\\
=:J_1+J_2+J_3
\end{aligned}
\end{equation}
where $J_1$, $J_2$, $J_3$ denote the integrations corresponding to
$\Delta_qR(B,v)$, $\Delta_qT_vB$, $\Delta_qT_Bv$.

\vspace{0.5cm}

Step 2. We have
\begin{equation}\label{4.4}
\sum_q|J_1| \leq C\|u\|_{L^{2,\infty}(0,1;L^\infty(\mathbb
R^3))}\|v\|_Q^2.
\end{equation}

\vspace{0.5cm}

As in the proof of Lemma \ref{2.1}, we denote
\begin{eqnarray*}
&& h(t)=\|\Delta_qv(t)\|_{L^2(\mathbb R^3)}\|\Delta_jv(t)\|_{L^2(\mathbb R^3)}\\
&& M=\sup_{0\leq t< 1}h(t)
\end{eqnarray*}
and define
$$
E_k=\{t\in(0,1):\quad 2^{-k}<M^{-1}h(t)\leq 2^{-{k-1}}\}
$$
where $M$ is bounded and $E_k$ are Lebesgue measurable because
$\|v\|_Q$ is bounded.

As in (\ref{2.7}) we have
\begin{equation}\label{4.5}
\begin{aligned}
 \int_0^1\|u(t)\|_{L^\infty(\mathbb R^3)}h(t)dt\\
 \leq C\|u\|_{L^{2,\infty}(0,1;L^\infty(\mathbb
R^3))}\sum_{k\geq-1}\frac1{|E_k|^{1/2}}\int_{E_k}h(t)dt\\
 \leq C\|u\|_{L^{2,\infty}(0,1;L^\infty(\mathbb
R^3))}\sum_{k\geq-1}(\sup_{E_k}h(t))^{1/2}(\int_{E_k}h(t))^{1/2}\\
 \leq  C\|u\|_{L^{2,\infty}(0,1;L^\infty(\mathbb
R^3))}M^{1/2}(\int_0^1h(t)dt)^{1/2}.
\end{aligned}\end{equation}
So
\begin{eqnarray*}
&& |J_1|\leq 2|\int_0^1dt\int_{\mathbb R^3}\Delta_q(\sum_{j\geq
q-3}\sum_{k=j-1}^{
j+1}\Delta_ku\Delta_jv)\cdot\nabla\Delta_qv(x,t)dx|\\
&& \leq C2^q\sum_{j\geq q-3}\int_0^1\|u(t)\|_{L^\infty(\mathbb
R^3)}\|\Delta_qv(t)\|_{L^2(\mathbb R^3)}\|\Delta_jv(t)\|_{L^2(\mathbb R^3)}dt\\
&& \leq C2^q\|u\|_{L^{2,\infty}(0,1;L^\infty(\mathbb
R^3))}\{\sup_{0\leq t<1}\|\Delta_qv(t)\|_{L^2(\mathbb R^3)}(\int_0^1\|\Delta_qv(t)\|_{L^2(\mathbb R^3)}^2dt)^{1/2}\}^{1/2}\\
&& \times\sum_{j\geq q-3}\{\sup_{0\leq
t<1}\|\Delta_jv(t)\|_{L^2(\mathbb
R^3)}(\int_0^1\|\Delta_jv(t)\|_{L^2(\mathbb R^3)}^2dt)^{1/2}\}^{1/2}
\end{eqnarray*}
by using (\ref{4.5}) and H\"{o}lder inequality, and
\begin{eqnarray*}
&& \sum_{q\geq-1}|J_1|\\
&& \leq C\|u\|_{L^{2,\infty}(0,1;L^\infty(\mathbb
R^3))}\left(\sum_{q\geq-1}2^q\sup_{0\leq t<1}\|\Delta_qv(t)\|_{L^2(\mathbb R^3)}(\int_0^1\|\Delta_qv(t)\|_{L^2(\mathbb R^3)}^2dt)^{1/2}\right)^{1/2}\\
&& \times\left(\sum_{q\geq-1}2^{q-3}(\sum_{j\geq q-3}\{\sup_{0\leq t<1}\|\Delta_jv(t)\|_{L^2(\mathbb R^3)}(\int_0^1\|\Delta_jv(t)\|_{L^2(\mathbb R^3)}^2dt)^{1/2}\}^{1/2})^2\right)^{1/2}\\
&& \leq C\|u\|_{L^{2,\infty}(0,1;L^\infty(\mathbb
R^3))}\left(\sum_{q\geq-1}\sup_{0\leq t<1}\|\Delta_qv(t)\|_{L^2(\mathbb R^3)}(\int_0^1\|\nabla\Delta_qv(t)\|_{L^2(\mathbb R^3)}^2dt)^{1/2}\right)^{1/2}\\
&& \times\left(\sum_{q\geq-1}(\sum_{j\geq q-3}2^{\frac{q-3-j}{2}}\{\sup_{0\leq t<1}\|\Delta_jv(t)\|_{L^2(\mathbb R^3)}(\int_0^1\|\nabla\Delta_jv(t)\|_{L^2(\mathbb R^3)}^2dt)^{1/2}\}^{1/2})^2\right)^{1/2}\\
&& +C\|u\|_{L^{2,\infty}(0,1;L^\infty(\mathbb R^3))}\sup_{0\leq t<1}\|\Delta_{-1}v(t)\|_{L^2(\mathbb R^3)}^2     \\
&& \leq C\|u\|_{L^{2,\infty}(0,1;L^\infty(\mathbb R^3))}\|v\|_Q^2,
\end{eqnarray*}
where the Hardy-Young inequality
$$
(\sum_{q\geq-1}(\sum_{j\geq q-2}2^{\frac{q-2-j}2}a_j)^2)^{1/2}
\leq C(\sum_{q\geq-1}a_q^2)^{1/2},\quad (a_j\geq 0),
$$
is used in the last step.

\vspace{0.5cm}

Step 3. We have
\begin{eqnarray*}
\sum_q |J_2| && =\sum_q |\int_0^1\int_{\mathbb
R^3}\Delta_q(\sum_{q-2\leq j\leq q+4}\Delta_ju S_{j-1}v)\cdot
\nabla\Delta_qv dxdt|\\
&& \leq C\sum_{-1\leq q\leq 1}\sum_{q-2\leq j\leq q+4}\int_0^1\|\Delta_ju(t)\|_{L^\infty(\mathbb R^3)}\|S_{j-1}v(t)\|_{L^2(\mathbb R^3)}\|\Delta_q\nabla v(t)\|_{L^2(\mathbb R^3)}dt\\
&& +C\sum_{q\geq 2} \sum_{q-2\leq j\leq q+4}\int_0^1\|\Delta_ju(t)\|_{L^2(\mathbb R^3)}\|S_{j-1}v(t)\|_{L^\infty(\mathbb R^3)}\|\Delta_q\nabla v(t)\|_{L^2(\mathbb R^3)}dt\\
&& \leq C\sum_{q\geq-1} \sum_{q-2\leq j\leq q+4}\int_0^1
\|\Delta_jv(t)\|_{L^2(\mathbb R^3)}\|S_{j-1}u(t)\|_{L^\infty(\mathbb
R^3)}\|\Delta_q\nabla v(t)\|_{L^2(\mathbb R^3)}dt
\end{eqnarray*}
because (see \cite{[C]})
\begin{eqnarray*}
&& \|S_{j-1}v(t)\|_{L^\infty(\mathbb R^3)}\leq C\|S_{j-1}\nabla u(t)\|_{L^\infty(\mathbb R^3)}\leq C2^j\|S_{j-1}u(t)\|_{L^\infty(\mathbb R^3)},\\
&& \|\Delta_ju(t)\|_{L^2(\mathbb R^3)}\leq C2^{-j}\|\Delta_j \nabla
u(t)\|_{L^2(\mathbb R^3)}\leq C2^{-j}\|\Delta_jv(t)\|_{L^2(\mathbb
R^3)}.
\end{eqnarray*}
So by using Lemma \ref{2.1}, we have
\begin{equation}\label{4.6}
\sum_q |J_2|\leq C\|u\|_{L^{2,\infty}(0,1;L^\infty(\mathbb
R^3))}\|v\|_Q^2.
\end{equation}

\vspace{0.5cm}

Step 4. Notice that
\begin{equation}\label{4.7}
\begin{aligned}
\sum_q|J_3|  =2\sum_q|\int_0^1\int_{\mathbb R^3}\Delta_q(\sum_{q-2\leq j\leq q+4}\Delta_jvS_{j-1}u)\cdot\nabla\Delta_qvdxdt|\\
\leq C\sum_q\sum_{q-2\leq j\leq q+4}\int_0^1\|S_{j-1}u(t)\|_{L^\infty(\mathbb R^3)}\|\Delta_jv(t)\|_{L^2(\mathbb R^3)}\|\Delta_q\nabla v(t)\|_{L^2(\mathbb R^3)}dt\\
\leq C\|u\|_{L^{2,\infty}(0,1;L^\infty(\mathbb R^3))}\|v\|_Q^2
\end{aligned}\end{equation}
by using Lemma \ref{2.1} again.

So from (\ref{4.3}), (\ref{4.4}), (\ref{4.6}) and (\ref{4.7}) we
have
$$
\|v\|_Q^2\leq C\|u\|_{L^{2,\infty}(0,1;L^\infty(\mathbb
R^3))}\|v\|_Q^2+\|v_0\|_{L^2(\mathbb R^3)}^2.
$$
Take $\epsilon=1/(2C)$. If $\|u\|_{L^{2,\infty}(0,1;L^\infty(\mathbb
R^3))}<\epsilon$, we have $\|v\|_Q^2\leq 2\|v_0\|_{L^2(\mathbb
R^3)}^2$.\qed

\end{document}